\documentclass[reqno,11pt]{amsart}
\usepackage[T1]{fontenc}
\usepackage{lmodern}
\usepackage[a4paper,margin=1in]{geometry}
\usepackage{amsmath,amssymb,amsthm,mathtools}
\usepackage{microtype}
\usepackage{xcolor}
\usepackage[hidelinks]{hyperref}
\usepackage{enumitem}
\setlist[enumerate]{itemsep=3pt,topsep=5pt}
\numberwithin{equation}{section}
\newtheorem{theorem}{Theorem}[section]
\newtheorem{lemma}[theorem]{Lemma}
\newtheorem{proposition}[theorem]{Proposition}
\newtheorem{corollary}[theorem]{Corollary}
\theoremstyle{remark}

\newcommand{\R}{\mathbb R}
\newcommand{\E}{\mathbb E}
\newcommand{\Prob}{\mathbb P}

\newcommand{\tr}{\operatorname{tr}}
\newcommand{\rank}{\operatorname{rank}}
\newcommand{\Ran}{\operatorname{Ran}}
\newcommand{\GDAN}{\mathrm{GDAN}}
\newcommand{\dto}{\xrightarrow{\mathrm d}}
\newcommand{\pto}{\xrightarrow{\mathrm P}}
\newcommand{\norm}[1]{\lVert #1\rVert}
\newcommand{\ip}[2]{\langle #1,#2\rangle}
\allowdisplaybreaks[2]
\hypersetup{pdftitle={Hotelling's T-squared and the generalized domain of attraction of the normal law},
  pdfsubject={Chi-square characterization and radial moment inequalities; version 8}}
\title[Hotelling's statistic and normal attraction]
{Hotelling's $T^2$ and the generalized domain of attraction of the normal law}
\date{}

\author{Wang Shaochen}
\address{School of Mathematics, South China University of Technology, Guangzhou 510640, China}
\email{mascwang@scut.edu.cn}

\author{Zhou Wang}
\address{Department of Statistics and Data Science, National University of Singapore, Singapore 119077, Singapore}
\email{stazw@nus.edu.sg}

\subjclass[2020]{Primary 60F05, 62H10; Secondary 62E20, 62H05}
\keywords{Hotelling's $T^2$; multivariate Student statistic; self-normalized sums;
  generalized domain of attraction of the normal law; quadratic negligibility.}

\begin{document}
\begin{abstract}
Let $X_1,X_2,\ldots$ be independent copies of a full random vector in
$\mathbb R^p$, where $p$ is fixed. We prove that Hotelling's $T^2$-statistic
at a prescribed location converges to $\chi_p^2$ if and only if the
underlying distribution belongs to the generalized domain of attraction
of the normal law and has that location as its mean. Our proof of
necessity combines moment bounds with matrix comparisons for
self-normalized sums. It yields a sharp moment inequality for the
squared norm of the uncentered self-normalized sum. Convergence of
the first two moments of this squared norm to the chi-square values
 implies generalized normal attraction and identifies the mean.
The necessity part of the multivariate Student normal-limit theorem
follows as a corollary. Its sufficiency is the established result of
Sepanski and of Vu, Maller and Klass.
\end{abstract}
\maketitle

\section{Introduction}

Hotelling's $T^2$-statistic is used to test a specified multivariate
mean and to construct confidence ellipsoids. Its chi-square calibration
is familiar when the dimension is fixed and the population covariance
matrix is finite and nonsingular. Generalized normal attraction extends
this calibration to some distributions with infinite second moments
\cite{Sepanski94}. The converse question asks which distributions can
produce the chi-square limit in the first place.

The one-dimensional characterization was proved by Gin\'e, G\"otze
and Mason \cite{GGM}, following work of Logan, Mallows, Rice and Shepp
\cite{LMRS} and the symmetric case of Griffin and Mason
\cite{GriffinMason}. O'Brien \cite{OBrien} characterized normal
attraction by the negligibility of the largest squared observation
relative to the sum of squares. For vectors, Hahn and Klass
\cite{HahnKlass} developed matrix normalization. Maller \cite{Maller}
and Sepanski \cite{SepanskiProb} gave probabilistic characterizations
of generalized normal attraction. Sepanski \cite{Sepanski94} proved
the Hotelling chi-square limit under generalized normal attraction
and the Student normal limit under ordinary normal attraction.
The latter result was extended to generalized normal attraction by
Sepanski \cite{Sepanski96} and by Vu, Maller and Klass \cite{VMK}.

Gin\'e and G\"otze \cite{GG} proved multivariate necessity for symmetric
vectors. They also obtained uniform radial moment bounds for nonsymmetric
laws assigning zero mass to every hyperplane through the origin.
Martsynyuk \cite[Theorem~3]{Martsynyuk} obtained converses when the
sample correlation matrix has a deterministic positive definite limit.
The symmetric-square-root version requires the identity limit.
For independent summands, Shao \cite{Shao18} obtained sufficient
conditions and proved necessity under maximal-term negligibility.
Li et al.~\cite{LLSSX} gave a characterization for independent triangular
arrays. Li, Pan and Zhou \cite{LPZ} study Hotelling $F$-approximations
with independent identically distributed entries and varying dimension.
Their symmetric characterization includes a fixed-dimensional
chi-square limit. They also give nonsymmetric sufficient conditions
when the dimension grows proportionally to the sample size.

In this paper, we establish the necessity of generalized normal attraction for the
chi-square limit of Hotelling's statistic under an arbitrary full
distribution in a fixed dimension. Combined with the sufficiency
theorem proved by
Sepanski~\cite{Sepanski96} and  Vu, Maller and Klass~\cite{VMK},
this completes the characterization. We also
solve the general multivariate Student converse problem discussed by
Giné and Götze \cite{GG}. Our proof addresses the mixed products in
the fourth radial moment identified in the
closing discussion of \cite{GG}. It also gives a sharp moment inequality
and a criterion based on the first two moments of the squared norm of
the uncentered self-normalized sum.

The paper is organized as follows. Section~\ref{sec:main} states the Hotelling and Student characterizations
and outlines the proof. Section~\ref{sec:preliminaries} reduces the
problem to an uncentered quadratic form, states the radial moment
criterion, and proves the moment bounds. Section~\ref{sec:radial-proof}
proves the criterion and completes the proofs of the main results.

\section{Main results}\label{sec:main}

Let $X,X_1,X_2,\ldots$ be independent and identically distributed (iid)
random vectors in $\R^p$, where $p\ge1$ is fixed. We assume the law of $X$ is full
but its coordinates may be dependent. Recall a law is \emph{full} if it is not supported by
an affine hyperplane. We write $\|\cdot\|$ for the Euclidean norm and
$\langle\cdot,\cdot\rangle$ for its inner product. The notation $A^\top$
denotes the transpose of a vector or matrix. We use $\dto$ and $\pto$ for
convergence in distribution and in probability, respectively. Write
\begin{equation}\label{eq:sample}
 \begin{gathered}
 S_n=\sum_{i=1}^nX_i,\qquad \bar X_n=S_n/n,\\
 \bar C_n=\sum_{i=1}^n(X_i-\bar X_n)(X_i-\bar X_n)^\top,
 \qquad \widehat\Sigma_n=\frac{\bar C_n}{n-1},\qquad n\ge2.
 \end{gathered}
\end{equation}
For a positive semidefinite matrix $A$, let $A^\dagger$ denote its
Moore--Penrose inverse and let $A^{\dagger/2}=(A^\dagger)^{1/2}$ be the
symmetric positive semidefinite square root. At a prescribed location
$\theta\in\R^p$, define Hotelling's statistic by
\begin{equation}\label{eq:hotelling}
 H_n(\theta)=n(\bar X_n-\theta)^\top\widehat\Sigma_n^\dagger
                       (\bar X_n-\theta).
\end{equation}
{The fullness} assumption of the law of $X$ makes $\bar C_n$ positive definite with probability tending
to one exponentially fast; see Lemma~\ref{lem:singular}. Thus the
convention on singular samples does not affect any weak limit below.

We say that $X$ belongs to the \emph{generalized domain of attraction of
the normal law} (GDAN), written $X\in\GDAN$, if there are deterministic
nonsingular matrices $A_n$ and vectors $b_n$ such that
\[
 A_n(S_n-b_n)\dto N_p(0,I_p).
\]
Here $I_p$ is the $p\times p$ identity matrix, and $N_p(0,I_p)$ is the
standard normal law on $\R^p$. In the ordinary normal domain of
attraction, the normalization is scalar, with a fixed nonsingular
limiting covariance matrix. We write $\chi_p^2$ for the chi-square law
with $p$ degrees of freedom.

Our main result below characterizes the chi-square limit at a prescribed location.

\begin{theorem}\label{thm:hotelling}
Suppose that the law of $X$ is full and let $\theta\in\R^p$. Then
\begin{equation}\label{eq:main-theorem}
 H_n(\theta)\dto\chi_p^2
 \quad\Longleftrightarrow\quad
 X\in\GDAN\ \text{and}\ \E X=\theta.
\end{equation}
\end{theorem}

The existence of the mean is part of the conclusion in the necessary
direction. Thus \eqref{eq:main-theorem} characterizes the full
distributions for which chi-square calibration is asymptotically valid
at every significance level.

To state the Student consequence, define
\begin{equation}\label{eq:student}
 \begin{gathered}
 C_n(\theta)=\sum_{i=1}^n(X_i-\theta)(X_i-\theta)^\top,\\
 T_n(\theta)=C_n(\theta)^{\dagger/2}(S_n-n\theta),\qquad
 \tau_n(\theta)=\bar C_n^{\dagger/2}(S_n-n\theta).
 \end{gathered}
\end{equation}
The usual normalization $\sqrt n\,\widehat\Sigma_n^{\dagger/2}
(\bar X_n-\theta)$ equals
$\sqrt{(n-1)/n}\,\tau_n(\theta)$, so it has the same weak limits.
For every sample,
\begin{equation}\label{eq:student-hotelling}
 H_n(\theta)=\frac{n-1}{n}\norm{\tau_n(\theta)}^2.
\end{equation}

\begin{corollary}\label{cor:student}
Suppose that the law of $X$ is full. For any fixed $\theta\in\R^p$,
the following conditions are equivalent:
\begin{enumerate}[label=\textup{(\roman*)}]
 \item $T_n(\theta)\dto N_p(0,I_p)$;
 \item $\tau_n(\theta)\dto N_p(0,I_p)$;
 \item $H_n(\theta)\dto\chi_p^2$;
 \item $X\in\GDAN$ and $\E X=\theta$;
 \item $\E\norm{T_n(\theta)}^2\to p$ and
       $\E\norm{T_n(\theta)}^4\to p^2+2p$.
\end{enumerate}
\end{corollary}

In particular, the necessity assertion for the centered Student vector
follows from Theorem~\ref{thm:hotelling} and
\eqref{eq:student-hotelling} by the continuous mapping theorem. The
passage from a radial limit to a vector limit uses the distributional
conclusion in Theorem~\ref{thm:hotelling} and the established Student
sufficiency theorem.

The proofs of Theorem~\ref{thm:hotelling} and
Corollary~\ref{cor:student} are given in
Section~\ref{subsec:completion}. We outline the proof of necessity here.
We first translate the prescribed location to the origin.
The chi-square limit is then transferred to the squared norm of the
uncentered self-normalized sum. An independent-block argument gives
uniform moment bounds from tightness. The first two moments therefore
converge to the chi-square values. Proposition~\ref{prop:radial} turns
these limits into quadratic negligibility. Maller's characterization \cite{Maller}
then gives generalized normal attraction, and integrability identifies
the mean.

The main step is the moment inequality developed in Proposition~\ref{prop:radial}.
An exact projection identity separates the individual quadratic
contributions from the mixed products. Conditional on the remaining
observations, a matrix-inverse kernel has an expectation equal to a
squared Hilbert-space norm. Differentiation gives the bound for products
involving three distinct observations. A rank-one perturbation of the
remaining second-moment matrix gives a one-sided bound for the other
mixed terms. Regularization extends the kernel argument to singular
samples. Section~\ref{subsec:identity} gives the projection identity and
the two estimates used to prove the inequality.

\section{Preliminary facts and moment bounds}\label{sec:preliminaries}

Translation preserves fullness and membership in $\GDAN$. We work at
$\theta=0$ until the final proofs. Thus
\begin{equation}\label{eq:basic}
 C_n=\sum_{i=1}^nX_iX_i^\top,\qquad
 T_n=C_n^{\dagger/2}S_n,\qquad
 Q_n=\norm{T_n}^2=S_n^\top C_n^\dagger S_n.
\end{equation}
The law is assumed full unless stated otherwise.

We write $a\preceq b$ when $a\le Cb$, for a finite constant independent
of $n$ and of any perturbation parameter. A subscript indicates other
allowed dependencies. Constants denoted by $c,C>0$ may change from line
to line. They may depend on the fixed dimension and the law of $X$.
For symmetric matrices, $A\le B$ means that $B-A$ is positive
semidefinite, and $A>0$ means positive definite. We write $\rank A$,
$\Ran A$ and $\tr A$ for the rank, range and trace of $A$. All sums over
distinct indices are ordered sums.

\subsection{Projection identities and radial reduction}\label{subsec:projection}

We first control singular samples. Fullness gives the following exponential bound.

\begin{lemma}\label{lem:singular}
If the law of $X$ is full, there are constants $c>0$ and $K<\infty$ such
that
\[
 \Prob(C_n\text{ is singular})+
 \Prob(\bar C_n\text{ is singular})\le K e^{-cn}.
\]
The same conclusion holds after deleting a fixed number of prescribed observations.
\end{lemma}

\begin{proof}
The support of a full law contains $p+1$ affinely independent points.
Choose open neighborhoods of these points so small that every choice of
one point from each neighborhood is affinely independent. Each
neighborhood has positive probability. Hence a block of $p+1$
observations has affine span $\R^p$ with some probability $\delta>0$.
If any one of the disjoint blocks has this property, the whole sample
has affine span $\R^p$. Its centered covariance matrix and its
uncentered second-moment matrix are then positive definite. The
probability that all $\lfloor n/(p+1)\rfloor$ blocks fail is at most
$(1-\delta)^{\lfloor n/(p+1)\rfloor}$. The assertion after deleting a
fixed number of observations follows in the same way.
\end{proof}

Define
\begin{equation}\label{eq:projection-notation}
 a_i=C_n^{\dagger/2}X_i,\qquad
 p_{ij}=a_i^\top a_j,\qquad P=(p_{ij})_{i,j\le n},\qquad
 \rho=\rank C_n.
\end{equation}
The dependence on $n$ will usually be suppressed. If
$\mathsf X_n=(X_1,\ldots,X_n)$ is the $p\times n$ data matrix, then
$P=\mathsf X_n^\top C_n^\dagger\mathsf X_n$ and
$P^2=P$, because $C_n^\dagger C_nC_n^\dagger=C_n^\dagger$.
Thus $P$ is the orthogonal projection onto the row space of $\mathsf X_n$, and
\begin{equation}\label{eq:projection}
 P^2=P=P^\top,\qquad
 0\le h_i:=p_{ii}\le1,\qquad
 \sum_i h_i=\sum_{i,j}p_{ij}^2=\rho\le p.
\end{equation}
We also write $\rho_n=\rho$ and $h_{n,i}=h_i$.
Writing
\begin{equation}\label{eq:q-v}
 t=\sum_i a_i,\qquad Q=\norm{t}^2=\norm{T_n}^2,\qquad
 v_i=\sum_jp_{ij},
\end{equation}
we also have
\begin{equation}\label{eq:projection-extra}
 \sum_i v_i=\sum_i v_i^2=Q,\qquad
 0\le Q\le n,\qquad
 |p_{ij}|\le\sqrt{h_i h_j}\le1.
\end{equation}
Let $v=(v_1,\ldots,v_n)^\top$ and $\mathbf 1_n=(1,\ldots,1)^\top$.
Then $v=P\mathbf 1_n$ and $Q=\mathbf 1_n^\top P\mathbf 1_n$.
The identities follow from $P^2=P$ and hold on singular samples as well. Define
\[
 D_n=\sum_{i=1}^n h_{n,i}^2.
\]
Then $0\le D_n\le p$, and
\begin{equation}\label{eq:leverage-comparison}
 \bigl(\max_{i\le n}h_{n,i}\bigr)^2
 \le D_n\le p\max_{i\le n}h_{n,i}.
\end{equation}

The next lemma relates Hotelling's statistic to the uncentered quadratic form.

\begin{lemma}\label{lem:radial-reduction}
On the event $\{\bar C_n>0\}$,
\begin{equation}\label{eq:radial-transform}
 Q_n=\frac{nH_n(0)}{n-1+H_n(0)},\qquad
 H_n(0)=\frac{(n-1)Q_n}{n-Q_n},\qquad Q_n<n.
\end{equation}
For a full law, tightness of either $(H_n(0))$ or $(Q_n)$ implies
$H_n(0)-Q_n\pto0$. In particular,
\begin{equation}\label{eq:radial-equivalence}
 H_n(0)\dto\chi_p^2
 \quad\Longleftrightarrow\quad Q_n\dto\chi_p^2.
\end{equation}
\end{lemma}

\begin{proof}
On this event let $U=S_n^\top\bar C_n^{-1}S_n$. Since
$C_n=\bar C_n+S_nS_n^\top/n$, the Sherman--Morrison formula gives
\[
 Q_n=\frac{U}{1+U/n},\qquad H_n(0)=\frac{n-1}{n}U.
\]
These formulas imply \eqref{eq:radial-transform}. The two rational maps
tend to the identity uniformly on bounded intervals. Thus tightness of
either statistic implies tightness of both. On this event,
\[
 H_n(0)-Q_n=\frac{Q_n(Q_n-1)}{n-Q_n}.
\]
Tightness and Lemma~\ref{lem:singular} now give convergence to zero in
probability.
\end{proof}

For a full law, Maller's quadratic-negligibility characterization
\cite{Maller}, as stated in \cite[Theorem~1.2(iv)]{GG}, gives
\begin{equation}\label{eq:maller}
 X\in\GDAN
 \quad\Longleftrightarrow\quad
 \max_{1\le i\le n}X_i^\top C_n^\dagger X_i\pto0.
\end{equation}
The pseudoinverse version is equivalent by Lemma~\ref{lem:singular}.
We also use the fact that $X\in\GDAN$ implies $\E\norm{X}<\infty$;
see \cite[Theorem~1.2(i)]{GG}. The following proposition gives a moment
criterion for quadratic negligibility.

\begin{proposition}\label{prop:radial}
Suppose that the law of $X$ is full and
\begin{equation}\label{eq:first-radial-moment}
 \E Q_n\longrightarrow p.
\end{equation}
Then
\begin{equation}\label{eq:radial-inequality}
 \limsup_{n\to\infty}\bigl\{\E Q_n^2+2\E D_n\bigr\}
 \le p^2+2p.
\end{equation}
If, in addition,
\begin{equation}\label{eq:radial-lower}
 \liminf_{n\to\infty}\E Q_n^2\ge p^2+2p,
\end{equation}
then
\begin{equation}\label{eq:D-small}
 \E D_n\longrightarrow0,
 \qquad X\in\GDAN,
 \qquad \E X=0.
\end{equation}
\end{proposition}

In particular, the two limits
\begin{equation}\label{eq:radial-moments}
 \E Q_n\longrightarrow p,\qquad
 \E Q_n^2\longrightarrow p^2+2p
\end{equation}
imply \eqref{eq:D-small}. The proposition also identifies $p^2+2p$ as
the largest possible asymptotic second radial moment under
\eqref{eq:first-radial-moment}. This value is attained when
$X\in\GDAN$ and $\E X=0$.

Proposition~\ref{prop:radial} is proved in
Section~\ref{subsec:completion}, using the estimates of
Sections~\ref{subsec:kernel} and \ref{subsec:perturbation}.
We next establish the moment bounds needed in that proof.
All moments of $Q_n$ exist for every $n$, since $0\le Q_n\le n$.

\subsection{Uniform moments from tightness}\label{subsec:moments}

We first prove a block estimate. It is useful to allow diagonal factors,
as well as factors with two distinct indices. We compare the common
sample denominator with separate block denominators by the following
inverse-order inequality.

\begin{lemma}\label{lem:inverse-order}
Let $A,B$ be positive semidefinite matrices with $A\le B$. If
$x\in\Ran A$, then
\[
 x^\top B^\dagger x\le x^\top A^\dagger x.
\]
\end{lemma}

\begin{proof}
The matrix order implies $\Ran A\subseteq\Ran B$. For any positive
semidefinite matrix $H$ and $x\in\Ran H$,
\[
 x^\top H^\dagger x
 =\sup_{z\in\R^p}\{2x^\top z-z^\top Hz\}.
\]
Apply this identity to $A$ and $B$.
\end{proof}

Applying this inequality to disjoint blocks gives the following estimate.
It extends \cite[Lemma~4.1]{GG} to singular blocks and applies to all iid laws.

\begin{lemma}\label{lem:block}
Let $r\ge1$ and consider a product
\[
 Z_n=\prod_{\ell=1}^r p_{\alpha_\ell\beta_\ell},
 \qquad \alpha_\ell,\beta_\ell\in\{1,\ldots,n\}.
\]
Suppose that exactly $m$ distinct indices occur in the list
\[
 \alpha_1,\beta_1,\ldots,\alpha_r,\beta_r.
\]
Let $s$ be the number of these distinct indices that occur
exactly once in this list. In particular, a factor $p_{ii}$
contributes two occurrences of $i$.
If $n\ge m$ and $b=\lfloor n/m\rfloor$, then
\begin{equation}\label{eq:block-bound}
 |\E Z_n|
 \le b^{-m}p^{m-s}\bigl(\E\norm{T_b}\bigr)^s.
\end{equation}
\end{lemma}

\begin{proof}
Relabel the $m$ distinct indices as $1,\ldots,m$, and let $r_\nu$ be
the number of occurrences of $\nu$. Choose disjoint sets
$J_1,\ldots,J_m\subseteq\{1,\ldots,n\}$, each of size $b$. By
exchangeability,
\begin{equation}\label{eq:block-averaging}
 b^m\E Z_n
 =\E\sum_{i_1\in J_1,\ldots,i_m\in J_m}
       \prod_{\ell=1}^r p_{i_{\alpha_\ell}i_{\beta_\ell}}.
\end{equation}
The matrix $C_n^\dagger$ is held fixed in this finite sum.

For each $\nu$ with $r_\nu=1$, linearity replaces $a_{i_\nu}$ by
\[
 A_\nu=\sum_{i\in J_\nu}a_i
       =C_n^{\dagger/2}S_{J_\nu},\qquad
 S_J=\sum_{i\in J}X_i.
\]
Apply the Cauchy--Schwarz inequality to every inner product after making
all these replacements. The absolute value of the sum in
\eqref{eq:block-averaging} is at most
\[
 \prod_{\nu:r_\nu=1}\norm{A_\nu}
 \prod_{\nu:r_\nu\ge2}\Big(\sum_{i\in J_\nu}h_i^{r_\nu/2}\Big).
\]
By \eqref{eq:projection}, each factor in the second product is at most
$p$. For the first product, let $C_J=\sum_{i\in J}X_iX_i^\top$.
Since $S_J\in\Ran C_J$ and $C_J\le C_n$,
Lemma~\ref{lem:inverse-order} gives
\[
 \norm{A_\nu}^2
 =S_{J_\nu}^\top C_n^\dagger S_{J_\nu}
 \le S_{J_\nu}^\top C_{J_\nu}^\dagger S_{J_\nu}.
\]
The square roots on the right are independent over the different
blocks and each has the law of $\norm{T_b}$. Taking expectations proves
\eqref{eq:block-bound}.
\end{proof}

The block estimate turns tightness into uniform bounds for all moments.

\begin{proposition}\label{prop:moments}
For any iid law on $\R^p$, if $(T_n)$ is tight, then
\begin{equation}\label{eq:all-moments}
 \sup_{n\ge1}\E\norm{T_n}^q<\infty\qquad\text{for every }q>0.
\end{equation}
Consequently, the family $(\norm{T_n}^q)_{n\ge1}$ is uniformly
integrable for every $q>0$.
\end{proposition}

\begin{proof}
Write $Y_n=\norm{T_n}$ and $\mu_n=\E Y_n$. These expectations are finite
because $Y_n^2\le n$. Applying Lemma~\ref{lem:block} to $p_{12}$, with
$b=\lfloor n/2\rfloor$, gives, for $n\ge2$,
\begin{equation}\label{eq:first-recurrence}
 \E Y_n^2
 =\E\rho+n(n-1)\E p_{12}
 \le p+\frac{n(n-1)}{b^2}\mu_b^2
 \le p+6\mu_b^2.
\end{equation}
Suppose that $(\mu_n)$ is unbounded. Choose a sequence $n_k$ at which
$\mu_{n_k}=\max_{j\le n_k}\mu_j$ and $\mu_{n_k}\to\infty$. The
Paley--Zygmund inequality and \eqref{eq:first-recurrence} yield, for all large $k$,
\[
 \Prob\left(Y_{n_k}\ge\frac{\mu_{n_k}}2\right)
 \ge\frac{\mu_{n_k}^2}{4\E Y_{n_k}^2}
 \ge\frac{1}{4(p+6)}.
\]
This contradicts tightness. Thus
\begin{equation}\label{eq:M}
 M:=\max\{1,\sup_n\E Y_n\}<\infty.
\end{equation}

Fix an integer $k\ge1$ and expand
\[
 Q_n^k=\left(\sum_{i,j=1}^n p_{ij}\right)^k.
\]
Group the terms according to the equality pattern of their $2k$ indices.
There are finitely many such patterns, depending only on $k$. A pattern
with $m$ distinct indices has at most $n^m$ assignments. If $s$ of its
indices occur once, Lemma~\ref{lem:block} bounds the absolute value of the
expected monomial by
\[
 \lfloor n/m\rfloor^{-m}p^{m-s}M^s.
\]
Since $n/\lfloor n/m\rfloor\le2m$ for $n\ge m$, multiplication by
$n^m$ gives a bound depending only on $k,p,M$ whenever $n\ge2k$. Summing over the patterns yields
$\sup_n\E Q_n^k<\infty$, noting that the finitely many smaller sample sizes satisfy
$Q_n\le n$. H\"older's inequality gives \eqref{eq:all-moments} for every
$q>0$. Applying the bound with an exponent larger than $q$ proves
uniform integrability.
\end{proof}

If $Q_n\dto\chi_p^2$, then $(T_n)$ is tight. Applying
Proposition~\ref{prop:moments} and uniform integrability gives
\eqref{eq:radial-moments}. Thus Lemma~\ref{lem:radial-reduction} reduces
the necessary part of Theorem~\ref{thm:hotelling} to
Proposition~\ref{prop:radial}.

\section{Quadratic negligibility from radial moments}\label{sec:radial-proof}

We prove Proposition~\ref{prop:radial} assuming fullness and
\eqref{eq:first-radial-moment}. The first moment bound makes $T_n$ tight.
Proposition~\ref{prop:moments} therefore gives
\begin{equation}\label{eq:derived-radial-bounds}
 \sup_{n\ge1}\E Q_n^r<\infty\qquad\text{for every }r>0.
\end{equation}
Since $\E\rho_n\to p$ by Lemma~\ref{lem:singular}, we have
\begin{equation}\label{eq:pair-small}
 \eta_n:=\E Q_n-\E\rho_n
       =n(n-1)\E p_{12}\longrightarrow0,\qquad n\ge2.
\end{equation}
The positive kernel in Lemma~\ref{lem:kernel} will also show that
$\eta_n\ge0$.

\subsection{Projection identity and proof outline}\label{subsec:identity}

We first give an identity for the individual quadratic contributions and
the mixed products. With the notation of
\eqref{eq:projection-notation}--\eqref{eq:q-v}, set
\begin{equation}\label{eq:r-q-W}
 \begin{gathered}
 r_{ij}=a_i^\top(t-a_i-a_j),\qquad
 q_{ij}=\norm{t-a_i-a_j}^2,\qquad i\ne j,\\
 W(\lambda)=\sum_{i\ne j}\frac{r_{ij}r_{ji}}{1+\lambda q_{ij}},
 \qquad \lambda\ge0,\\
 D=\sum_i h_i^2,\qquad
 B_* =\sum_{\substack{i,j,k\\\text{pairwise distinct}}}p_{ij}p_{ik}.
 \end{gathered}
\end{equation}
For the sample projection, $D=D_n$. We also write $W_n(\lambda)$ and
$B_{*,n}$.
Here $r_{ij}$ is the inner product of the normalized observation $i$
and the sum with observations $i,j$ removed. Expanding this sum gives
\[
 W(0)=B_*+\sum_{i,j,k,l}^{\ne}p_{ik}p_{jl}.
\]
Thus $W(0)$ contains the four-index contribution together with one copy
of the three-index sum.

\begin{lemma}\label{lem:identity}
For any orthogonal projection $P=(p_{ij})$ of rank $\rho$, the
quantities in \eqref{eq:r-q-W} satisfy
\begin{equation}\label{eq:main-identity}
 W(0)+3B_*=(Q-\rho)^2-2\rho+2D.
\end{equation}
\end{lemma}

\begin{proof}
For each fixed $i$, subtract the terms with $j=k$ from
$(\sum_{j\ne i}p_{ij})^2$. Since $\sum_jp_{ij}^2=h_i$, this gives
\begin{align}\label{eq:Bidentity}
 B_*&=\sum_i(v_i-h_i)^2-\sum_{i\ne j}p_{ij}^2\notag\\
 &=Q-\rho-2\sum_i h_iv_i+2D.
\end{align}
With $w_i=v_i-h_i$, expand
\[
 W(0)=\sum_{i\ne j}(w_i-p_{ij})(w_j-p_{ij}).
\]
The first part is $(\sum_iw_i)^2-\sum_iw_i^2$, the two mixed parts
sum to $-2\sum_iw_i^2$, and the last part is $\rho-D$. Therefore
\begin{align}\label{eq:Widentity}
 W(0)
 &=(Q-\rho)^2-3\sum_i(v_i-h_i)^2+\rho-D\notag\\
 &=(Q-\rho)^2-3Q+6\sum_i h_iv_i+\rho-4D.
\end{align}
Adding three times \eqref{eq:Bidentity} to \eqref{eq:Widentity} proves
\eqref{eq:main-identity}.
\end{proof}

We will prove that, for $n\ge3$ and $0<\lambda\le1/4$,
\begin{equation}\label{eq:two-key-estimates}
 |\E B_{*,n}|\le C\sqrt{\eta_n},\qquad
 \E W_n(0)\le\frac{\eta_n}{\lambda}+C\lambda.
\end{equation}
The constant is independent of $n$ and $\lambda$. Since $\eta_n\to0$,
these estimates give $\E B_{*,n}\to0$ and
$\limsup_n\E W_n(0)\le0$. The identity
\eqref{eq:main-identity} then yields the radial moment inequality.
Proposition~\ref{prop:triple} proves the first estimate in
Section~\ref{subsec:kernel}. Proposition~\ref{prop:W} proves the second
in Section~\ref{subsec:perturbation}.

\subsection{A positive kernel and the three-index term}\label{subsec:kernel}

For a positive definite matrix $B$, consider
\[
 K_{B,t}(x,y)=x^\top(tB+xx^\top+yy^\top)^{-1}y,\qquad t>0.
\]
The following lemma gives a Hilbert-space representation of this kernel.
For independent copies, its expectation is the squared norm of a mean
vector.

\begin{lemma}\label{lem:kernel}
For every $B>0$, there is a differentiable Hilbert-space map
$t\mapsto\Phi_t^B(x)$ such that
\begin{equation}\label{eq:feature-norm}
 K_{B,t}(x,y)=\ip{\Phi_t^B(x)}{\Phi_t^B(y)},\qquad
 \sup_x\norm{\Phi_t^B(x)}^2\le\frac12,\qquad
 \sup_x\norm{\partial_t\Phi_t^B(x)}^2\le\frac{1}{8t^2}.
\end{equation}
Consequently, for iid $X,Y$,
\begin{equation}\label{eq:kernel-positive}
 \E K_{B,t}(X,Y)=\norm{\E\Phi_t^B(X)}^2\ge0.
\end{equation}
For $C=B+xx^\top+yy^\top$, let $\kappa=x^\top C^{-1}y$,
$h_x=x^\top C^{-1}x$ and $h_y=y^\top C^{-1}y$. Then
\begin{equation}\label{eq:derivative-kernel}
 \ip{\partial_t\Phi_t^B(x)|_{t=1}}{\partial_t\Phi_t^B(y)|_{t=1}}
 =\kappa\bigl[(h_x-\tfrac12)(h_y-\tfrac12)+\kappa^2\bigr].
\end{equation}
For $B\ge0$, the kernel
\[
 K_{B,t}(x,y)=x^\top(tB+xx^\top+yy^\top)^\dagger y
\]
is also positive semidefinite, is bounded in absolute value by $1/2$,
and satisfies $\E K_{B,t}(X,Y)\ge0$ for iid $X,Y$.
\end{lemma}

\begin{proof}
Set $z=B^{-1/2}x$, $w=B^{-1/2}y$, $a=\norm{z}^2$, $b=\norm{w}^2$ and
$c=z^\top w$. Matrix inversion gives
\begin{equation}\label{eq:kernel-scalar}
 K_{B,t}(x,y)=\frac{tc}{(t+a)(t+b)-c^2}.
\end{equation}
Let $\mathcal H$ be the Hilbert direct sum of the real tensor spaces
$(\R^p)^{\otimes(2k+1)}$, $k\ge0$, with their usual tensor inner
products. We use the same inner-product and norm notation on
$\mathcal H$. For a positive integer $m$, $z^{\otimes m}$ denotes the
$m$-fold tensor product of $z$ with itself. Define
\begin{equation}\label{eq:feature}
 \Phi_t^B(x)=\bigoplus_{k=0}^\infty
       \frac{\sqrt t\,z^{\otimes(2k+1)}}{(t+a)^{k+1}}.
\end{equation}
Since $c^2\le ab<(t+a)(u+b)$ for $t,u>0$, summing a geometric series yields
\begin{equation}\label{eq:mixed-kernel}
 \ip{\Phi_t^B(x)}{\Phi_u^B(y)}
 =\frac{\sqrt{tu}\,c}{(t+a)(u+b)-c^2}.
\end{equation}
In particular, $\norm{\Phi_t^B(x)}^2=a/(t+2a)\le1/2$, and
\eqref{eq:kernel-scalar} gives the kernel representation. The feature map
is measurable and bounded in a separable Hilbert space. Its Bochner mean
is therefore defined, and independence gives \eqref{eq:kernel-positive}.

For fixed $x$, the series in \eqref{eq:feature} and its derivative
converge locally uniformly in the Hilbert-space norm for $t>0$. The
$k$th derivative component is
\[
 \left(\frac{1}{2t}-\frac{k+1}{t+a}\right)
 \frac{\sqrt t\,z^{\otimes(2k+1)}}{(t+a)^{k+1}}.
\]
Write $D_{t,u}=(t+a)(u+b)-c^2$ and let $F(t,u)$ denote the
right-hand side of \eqref{eq:mixed-kernel}. Differentiating once in each
variable gives
\[
 \partial_t\partial_u F
 =F\left[
 \left(\frac1{2t}-\frac{u+b}{D_{t,u}}\right)
 \left(\frac1{2u}-\frac{t+a}{D_{t,u}}\right)
 +\frac{c^2}{D_{t,u}^2}\right].
\]
At $t=u=1$, let $D_0=D_{1,1}$. The formulas
\[
 \kappa=\frac{c}{D_0},\qquad
 h_x=1-\frac{1+b}{D_0},\qquad
 h_y=1-\frac{1+a}{D_0}
\]
give \eqref{eq:derivative-kernel}. Taking $x=y$ and rescaling $t$,
we obtain
\[
 \norm{\partial_t\Phi_t^B(x)}^2
 =t^{-2}h_t\bigl[(h_t-\tfrac12)^2+h_t^2\bigr],
 \qquad h_t=\frac{a}{t+2a}\in[0,1/2].
\]
Since $(h-1/2)^2+h^2\le1/4$ on $[0,1/2]$, this is at most
$1/(8t^2)$.

Finally, let $B\ge0$ and replace it by $B+\varepsilon I_p$.
Both $x$ and $y$ belong to the range of $tB+xx^\top+yy^\top$, so
spectral decomposition gives
\[
 x^\top\{tB+t\varepsilon I_p+xx^\top+yy^\top\}^{-1}y
 \longrightarrow K_{B,t}(x,y).
\]
The positive definite case bounds the kernels by $1/2$. Their finite
kernel matrices are positive semidefinite, and this property is
preserved under passage to the limit. Bounded convergence gives the
expectation inequality.
\end{proof}

The first radial moment controls the squared norm of the conditional
mean in \eqref{eq:kernel-positive}. The block estimate controls the
squared norm of its derivative. Their inner product gives the
three-index estimate.

\begin{proposition}\label{prop:triple}
Under \eqref{eq:first-radial-moment}, $\eta_n\ge0$ and, for $n\ge3$,
\begin{equation}\label{eq:triple-small}
 |\E B_{*,n}|\le C\sqrt{\eta_n}.
\end{equation}
In particular, $\E(p_{12}p_{13})=o(n^{-3})$.
\end{proposition}

\begin{proof}
Let $B=\sum_{i=3}^nX_iX_i^\top$. The first two observations are iid
and independent of $B$. Conditioning on $B$ in Lemma~\ref{lem:kernel}
gives $\E p_{12}\ge0$, and hence $\eta_n\ge0$.

For $\varepsilon>0$, set
\[
 B_\varepsilon=B+\varepsilon I_p,\qquad
 C_\varepsilon=C_n+\varepsilon I_p,
 \qquad
 \kappa_\varepsilon=X_1^\top C_\varepsilon^{-1}X_2,
 \qquad h_{\varepsilon,i}=X_i^\top C_\varepsilon^{-1}X_i.
\]
For fixed $B$, let
$m_\varepsilon(t)=\E_X\Phi_t^{B_\varepsilon}(X)$ and
$f_\varepsilon(t)=\norm{m_\varepsilon(t)}^2$.
Here $\E_X$ integrates an independent copy $X$ while $B$ is held fixed;
outer expectations also integrate $B$. The uniform bounds in
\eqref{eq:feature-norm} justify differentiation under $\E_X$ locally in
$t$. Independence and \eqref{eq:derivative-kernel} give
\begin{align*}
 \alpha_{n,\varepsilon}
 &:=\E\norm{m_\varepsilon(1)}^2=\E\kappa_\varepsilon,\\
 \beta_{n,\varepsilon}
 &:=\E\norm{m'_\varepsilon(1)}^2\\
 &=\E\left[\kappa_\varepsilon
   \left\{(h_{\varepsilon,1}-\tfrac12)
          (h_{\varepsilon,2}-\tfrac12)
          +\kappa_\varepsilon^2\right\}\right].
\end{align*}
Since $f'_\varepsilon(1)=2\ip{m_\varepsilon(1)}{m'_\varepsilon(1)}$,
the Cauchy--Schwarz inequality gives
$|\E f'_\varepsilon(1)|\le
2\sqrt{\alpha_{n,\varepsilon}\beta_{n,\varepsilon}}$.
Differentiating the inverse kernel therefore yields
\begin{equation}\label{eq:regularized-derivative}
 \left|\E X_1^\top C_\varepsilon^{-1}B_\varepsilon
                         C_\varepsilon^{-1}X_2\right|
 \le2\sqrt{\alpha_{n,\varepsilon}\beta_{n,\varepsilon}}.
\end{equation}

For $C_{\varepsilon,t}=tB_\varepsilon+X_1X_1^\top+X_2X_2^\top$,
the order $B_\varepsilon\le t^{-1}C_{\varepsilon,t}$ bounds the
absolute value of the derivative of the kernel integrand by $t^{-1}$.
This also justifies differentiation under the expectation.

We now let $\varepsilon\downarrow0$ with $n$ fixed. All observations lie
in $\Ran C_n$, and $\Ran B\subseteq\Ran C_n$. Spectral decomposition
therefore gives $C_\varepsilon^{-1}X_i\to C_n^\dagger X_i$ and
\[
 X_1^\top C_\varepsilon^{-1}B_\varepsilon C_\varepsilon^{-1}X_2
 \longrightarrow\sum_{k=3}^n p_{1k}p_{2k}.
\]
By $B_\varepsilon\le C_\varepsilon$ and
$0\le h_{\varepsilon,i}\le1$, the absolute value of the expression on the
left is at most one. The integrands defining $\alpha_{n,\varepsilon}$ and
$\beta_{n,\varepsilon}$ are also bounded uniformly in $\varepsilon$. Bounded convergence in
\eqref{eq:regularized-derivative} gives
\begin{equation}\label{eq:fprime-small}
 (n-2)|\E(p_{13}p_{23})|
 \le2\sqrt{(\E p_{12})\beta_n},
\end{equation}
where
\begin{equation}\label{eq:beta}
 \beta_n:=\E\left[p_{12}
       \left\{(h_1-\tfrac12)(h_2-\tfrac12)+p_{12}^2\right\}\right]
 \ge0.
\end{equation}

The block estimate bounds the absolute value of each expected product
in the expansion of \eqref{eq:beta}. The products $p_{12}$, $p_{12}h_1$,
$p_{12}h_1h_2$ and $p_{12}^3$ involve two distinct indices, with
multiplicities $(1,1)$, $(3,1)$, $(3,3)$ and $(3,3)$, respectively.
Lemma~\ref{lem:block} and the bounded first moments of $\norm{T_n}$
thus give
\begin{equation}\label{eq:beta-bound}
 0\le\beta_n\le Cn^{-2}.
\end{equation}
Since $\E p_{12}=\eta_n/\{n(n-1)\}$, multiplying
\eqref{eq:fprime-small} by $n(n-1)$ proves
\[
 n(n-1)(n-2)|\E(p_{13}p_{23})|\le C\sqrt{\eta_n}.
\]
Relabeling the three indices and using exchangeability proves
\eqref{eq:triple-small}. The last assertion follows from $\eta_n\to0$.
\end{proof}

\subsection{Rank-one perturbation}\label{subsec:perturbation}

The perturbation $B\mapsto B+\lambda ss^\top$, where $s$ is the sum of
the remaining observations, gives
$\E W_n(\lambda)\le\eta_n/\lambda$.
To pass to $W_n(0)$, we use the following deterministic estimate. Its
constant does not depend on the sample size. The proof sums the signed
coefficients before taking absolute values.

\begin{lemma}\label{lem:remainder}
Let $a_1,\ldots,a_n\in\R^p$ have an orthogonal projection as their
Gram matrix. With $Q$ and $W(\lambda)$ as in
\eqref{eq:q-v} and \eqref{eq:r-q-W}, for $0<\lambda\le1/4$,
\begin{equation}\label{eq:remainder}
 |W(0)-W(\lambda)|\preceq_p\lambda(1+Q)^3.
\end{equation}
\end{lemma}

\begin{proof}
All the projection identities \eqref{eq:projection} and
\eqref{eq:projection-extra} apply, with $\rho=\rank P\le p$. Set
\[
 w_i=v_i-h_i,\qquad u_i=h_i-2v_i,\qquad
 L=\sum_iw_i=Q-\rho,\qquad f(x)=\frac{x}{1+\lambda x}.
\]
Write $w=(w_1,\ldots,w_n)^\top$ and $u=(u_1,\ldots,u_n)^\top$.
The following bounds will be used repeatedly:
\begin{equation}\label{eq:wu-bounds}
 \begin{gathered}
 |L|\le Q+p,\qquad \norm{w}^2\le2Q+2p,\qquad
 \norm{u}^2\le8Q+2p,\\
 \norm{v}^2=Q,\qquad \max_i|w_i|\le\sqrt Q+1,\qquad
 \sum_{i,j}p_{ij}^2=\rho\le p.
 \end{gathered}
\end{equation}
These follow from $\sum_i h_i^2\le p$ and $|v_i|\le\sqrt Q$.

Extend $r_{ij}$ and $q_{ij}$ to all $i,j$ by the formulas in
\eqref{eq:r-q-W}. Then
\begin{equation}\label{eq:rq-expand}
 r_{ij}=w_i-p_{ij},\qquad
 q_{ij}=Q+u_i+u_j+2p_{ij}.
\end{equation}
Also,
\[
 Q+u_i=\norm{t-a_i}^2\ge0,\qquad
 Q+u_i+u_j=q_{ij}-2p_{ij}\ge-2.
\]
Every value $Q+\eta u_i+\zeta u_j$, $0\le\eta,\zeta\le1$, is a
convex combination of the four corner values
$Q,Q+u_i,Q+u_j,Q+u_i+u_j$. It is therefore at least $-2$, even
when $u_i$ or $u_j$ is negative. On $[-2,\infty)$,
\[
 |f'(x)|\le4,\qquad |f''(x)|\le4,
\]
uniformly in $0<\lambda\le1/4$. A double application of the fundamental
theorem of calculus, followed by the change $2p_{ij}$ in the argument,
gives
\begin{equation}\label{eq:mixed-taylor}
 f(q_{ij})=f(Q+u_i)+f(Q+u_j)-f(Q)+e_{ij},\qquad
 |e_{ij}|\le4|u_i u_j|+8|p_{ij}|.
\end{equation}

We estimate the full sum
\[
 F=\sum_{i,j}(w_i-p_{ij})(w_j-p_{ij})f(q_{ij})
   =F_1-2F_2+F_3,
\]
where
\[
 F_1=\sum_{i,j}w_iw_jf(q_{ij}),\quad
 F_2=\sum_{i,j}w_i p_{ij}f(q_{ij}),\quad
 F_3=\sum_{i,j}p_{ij}^2f(q_{ij}).
\]
For $F_1$, let
$\Delta=\sum_iw_i\{f(Q+u_i)-f(Q)\}$. Equation~\eqref{eq:mixed-taylor}
gives the exact decomposition
\[
 F_1=f(Q)L^2+2L\Delta+\sum_{i,j}w_iw_je_{ij}.
\]
We sum the terms to obtain $L=Q-\rho$ before taking absolute values.
Since
$|\Delta|\le4\norm{w}\norm{u}$, equation~\eqref{eq:mixed-taylor} yields
\begin{align}\label{eq:F1}
 |F_1|
 &\le QL^2+8|L|\norm{w}\norm{u}
       +4\norm{w}^2\norm{u}^2
       +8\sqrt p\,\norm{w}^2\notag\\
 &\preceq_p(1+Q)^3.
\end{align}
Here we used
\[
 \sum_{i,j}|w_iw_jp_{ij}|
 \le\left(\sum_{i,j}p_{ij}^2\right)^{1/2}
      \left(\sum_{i,j}w_i^2w_j^2\right)^{1/2}
 \le\sqrt p\,\norm{w}^2.
\]

For $F_2$, both $q_{ij}$ and $Q+u_i$ are nonnegative. On that interval
$f$ is one-Lipschitz, so
\[
 |f(q_{ij})-f(Q+u_i)|\le |u_j|+2|p_{ij}|.
\]
Using $\sum_jp_{ij}=v_i$, the contribution with $f(Q+u_i)$ is
$\sum_iw_iv_i f(Q+u_i)$. Since
$0\le f(Q+u_i)\le(\sqrt Q+1)^2$, we obtain
\begin{align}\label{eq:F2}
 |F_2|
 &\le (\sqrt Q+1)^2\norm{w}\sqrt Q
       +\sqrt p\,\norm{w}\norm{u}
       +2p\max_i|w_i|\notag\\
 &\preceq_p(1+Q)^2.
\end{align}
The second term follows by the Cauchy--Schwarz inequality over the pair $(i,j)$.
The third uses $\sum_{i,j}p_{ij}^2\le p$.

Finally, $0\le f(q_{ij})\le q_{ij}\le(\sqrt Q+2)^2$, and hence
\begin{equation}\label{eq:F3}
 0\le F_3\le p(\sqrt Q+2)^2.
\end{equation}
Equations~\eqref{eq:F1}--\eqref{eq:F3} show that
$|F|\preceq_p(1+Q)^3$.

The diagonal contribution to $F$ is bounded by
\[
 \sum_i(v_i-2h_i)^2f(q_{ii})
 \le(2Q+8p)(\sqrt Q+2)^2\preceq_p(1+Q)^2.
\]
Thus the same cubic bound holds after restricting the sum to $i\ne j$.
By \eqref{eq:rq-expand},
\[
 W(0)-W(\lambda)
 =\lambda\sum_{i\ne j}(w_i-p_{ij})(w_j-p_{ij})f(q_{ij}),
\]
which proves \eqref{eq:remainder}.
\end{proof}

The preceding lemma allows the perturbation to be removed after taking
expectations. We first bound the perturbed quantity by conditioning on
the remaining observations and using kernel positivity.

\begin{proposition}\label{prop:W}
Under \eqref{eq:first-radial-moment}, for $n\ge2$ and
$0<\lambda\le1/4$,
\begin{equation}\label{eq:Wlimsup}
 \E W_n(0)\le\frac{\eta_n}{\lambda}+C\lambda.
\end{equation}
Consequently, $\limsup_n\E W_n(0)\le0$.
\end{proposition}

\begin{proof}
Fix $\lambda>0$ and let
\[
 B=\sum_{i=3}^nX_iX_i^\top,\qquad
 s=\sum_{i=3}^nX_i,\qquad A=B+\lambda ss^\top.
\]
For $n=2$, the two sums are zero.
Conditioning on $X_3,\ldots,X_n$ and applying the semidefinite part of
Lemma~\ref{lem:kernel} give
\[
 \E\left[X_1^\top(A+X_1X_1^\top+X_2X_2^\top)^\dagger X_2\right]
 \ge0.
\]
Since $s\in\Ran C_n$, the matrices $C_n$ and $C_n+\lambda ss^\top$
have the same range. The Sherman--Morrison formula on this range gives
\begin{equation}\label{eq:SM}
 X_1^\top(C_n+\lambda ss^\top)^\dagger X_2
 =p_{12}-\lambda\frac{r_{12}r_{21}}{1+\lambda q_{12}}.
\end{equation}
All terms are integrable, since $|p_{12}|\le1$ and
\[
 |r_{12}r_{21}|\le\sqrt{h_1h_2}\,q_{12}\le q_{12},\qquad
 \left|\frac{r_{12}r_{21}}{1+\lambda q_{12}}\right|\le\lambda^{-1}.
\]
Taking expectations in \eqref{eq:SM} and using exchangeability, we obtain
\begin{equation}\label{eq:Wlambda}
 \E W_n(\lambda)\le\frac{n(n-1)\E p_{12}}{\lambda}
                   =\frac{\eta_n}{\lambda},\qquad \lambda>0.
\end{equation}
For $0<\lambda\le1/4$, Lemma~\ref{lem:remainder} and
\eqref{eq:derived-radial-bounds} give
\[
 \E|W_n(0)-W_n(\lambda)|\le C\lambda\E(1+Q_n)^3\le C\lambda.
\]
Combining the two estimates proves \eqref{eq:Wlimsup}. First let
$n\to\infty$ with $\lambda$ fixed, and then let $\lambda\downarrow0$.
\end{proof}

\subsection{Completion of the proofs}\label{subsec:completion}

The two estimates now give the radial moment inequality. We then use it
to prove the Hotelling and Student characterizations.

\begin{proof}[Proof of Proposition~\ref{prop:radial}]
Taking expectations in \eqref{eq:main-identity}, we obtain
\begin{equation}\label{eq:closure}
 \E Q_n^2+2\E D_n
 =2\E(Q_n\rho_n)-\E\rho_n^2+2\E\rho_n
   +\E W_n(0)+3\E B_{*,n}.
\end{equation}
To bound the rank terms, let $R_n=p-\rho_n\ge0$.
Since $\E Q_n=p-\E R_n+\eta_n$,
\begin{align*}
 2\E(Q_n\rho_n)-\E\rho_n^2+2\E\rho_n
 &=p^2+2p+2p\eta_n
   -2\E(Q_nR_n)-\E R_n^2-2\E R_n\\
 &\le p^2+2p+2p\eta_n.
\end{align*}
Propositions~\ref{prop:triple} and \ref{prop:W} therefore imply, for
$n\ge3$ and $0<\lambda\le1/4$,
\begin{equation}\label{eq:radial-quantitative}
 \E Q_n^2+2\E D_n
 \le p^2+2p+2p\eta_n+C\sqrt{\eta_n}
          +\frac{\eta_n}{\lambda}+C\lambda.
\end{equation}
For $0<\eta_n\le1/16$, choose $\lambda=\sqrt{\eta_n}$. The excess
over $p^2+2p$ is then at most $C\sqrt{\eta_n}$. If $\eta_n=0$,
let $\lambda\downarrow0$ in \eqref{eq:radial-quantitative}. Since
$\eta_n\to0$, this proves \eqref{eq:radial-inequality}.

If \eqref{eq:radial-lower} also holds, then
\[
 0\le\limsup_{n\to\infty}2\E D_n
 \le p^2+2p-\liminf_{n\to\infty}\E Q_n^2\le0.
\]
Thus $\E D_n\to0$. For every $\varepsilon>0$,
\[
 \Prob\left(\max_{i\le n}h_{n,i}>\varepsilon\right)
 \le\varepsilon^{-2}\E D_n\longrightarrow0.
\]
Maller's characterization \eqref{eq:maller} gives $X\in\GDAN$.

It remains to identify the mean. Generalized normal attraction implies
$\E\norm{X}<\infty$. Integrability gives
$n\Prob(\norm{X}>\varepsilon n)\to0$, and hence
$\max_{i\le n}\norm{X_i}/n\pto0$. The law of large numbers now gives
\begin{equation}\label{eq:trace-small}
 \frac{\tr C_n}{n^2}
 \le\left(\max_{i\le n}\frac{\norm{X_i}}n\right)
       \left(\frac1n\sum_{i=1}^n\norm{X_i}\right)
 \pto0.
\end{equation}
Since $S_n\in\Ran C_n$, the spectral decomposition of $C_n$ gives
$\norm{S_n}^2\le(\tr C_n)Q_n$, also on singular samples. Tightness of
$Q_n$ and \eqref{eq:trace-small} imply $S_n/n\pto0$.
By the law of large numbers, this limit is $\E X$. Hence $\E X=0$.
\end{proof}

We now complete the proofs of the two characterizations. The radial
moment criterion gives necessity. The Student normal-limit theorem
gives sufficiency.

\begin{proof}[Proof of Theorem~\ref{thm:hotelling}]
First take $\theta=0$. If $H_n(0)\dto\chi_p^2$,
Lemma~\ref{lem:radial-reduction} gives $Q_n\dto\chi_p^2$.
Proposition~\ref{prop:moments} supplies \eqref{eq:radial-moments}, and
Proposition~\ref{prop:radial} yields $X\in\GDAN$ and $\E X=0$.
This proves necessity using only the radial limit.

Conversely, suppose that $X\in\GDAN$ and $\E X=0$. The established
Student sufficiency theorem \cite{Sepanski96,VMK}, in the form
\cite[Theorem~1.3]{GG}, gives
\[
 T_n\dto N_p(0,I_p),\qquad \tau_n(0)\dto N_p(0,I_p).
\]
Equation~\eqref{eq:student-hotelling} gives
$H_n(0)\dto\chi_p^2$.

For general $\theta$, apply the result to $X-\theta$. Translation
preserves fullness,  $\GDAN$, and $\bar C_n$.
The mean conclusion becomes $\E X=\theta$.
\end{proof}

It remains to connect the vector and moment conditions in
Corollary~\ref{cor:student}.

\begin{proof}[Proof of Corollary~\ref{cor:student}]
Theorem~\ref{thm:hotelling} proves (iii)$\Longleftrightarrow$(iv).
The Student sufficiency theorem used above, applied to $X-\theta$,
proves (iv)$\Longrightarrow$(i) and (ii).

If (ii) holds, \eqref{eq:student-hotelling} and the continuous mapping
theorem give (iii). If (i) holds, the same theorem gives
$\norm{T_n(\theta)}^2\dto\chi_p^2$, and
Lemma~\ref{lem:radial-reduction}, applied to the translated sample,
again gives (iii).

Condition (i), Proposition~\ref{prop:moments}, and uniform integrability
give (v). Conversely, Proposition~\ref{prop:radial}, applied to
$X-\theta$, gives (v)$\Longrightarrow$(iv). This proves all five
equivalences.
\end{proof}

\end{document}